\documentclass[12pt]{article}

\usepackage{amsfonts,amssymb,amsmath,color}

\newcommand\al{\alpha}
\newcommand\bt{\beta}
\newcommand\dl{\delta}
\newcommand\eps{\epsilon}
\newcommand\lm{\lambda}
\newcommand\sg{\sigma}
\newcommand\Gm{\Gamma}
\newcommand\Lm{\Lambda}
\newcommand\Sg{\Sigma}

\newcommand\F{\mathbb{F}}

\newcommand\cF{\mathcal{F}}
\newcommand\cJ{\mathcal{J}}
\newcommand\cM{\mathcal{M}}

\newcommand\ad{\mathrm{ad}}
\newcommand\Aut{\mathrm{Aut}}

\newcommand\la{\langle}
\newcommand\ra{\rangle}
\newcommand\lla{\langle\!\langle}
\newcommand\rra{\rangle\!\rangle}
\newcommand\llf{(\!(}
\newcommand\rrf{)\!)}

\newtheorem{theorem}{Theorem}[section]
\newtheorem{example}[theorem]{Example}
\newtheorem{proposition}[theorem]{Proposition}

\newtheorem{corollary}[theorem]{Corollary}
\newtheorem{definition}[theorem]{Definition}

\newcommand\pf{\noindent{\bf Proof.} }
\newcommand\qed{\hfill$\Box$}

\title{Radicals in flip subalgebras}
\author{B. G. Rodrigues\thanks{This work is based on the research supported by the National Research Foundation of South Africa (Grant Number CPRR23041894647)}\\\small Department of Mathematics and Applied Mathematics\\\small University of Pretoria\\\small Hatfield 0028, Pretoria, South Africa\\\small bernardo.rodrigues@up.ac.za\bigskip\\S.~Shpectorov \\\small{School of Mathematics}\\\small University of Birmingham\\\small Edgbaston, Birmingham B15 2TT, United Kingdom\\\small  s.shpectorov@bham.ac.uk\\
\small and\\
\small Department of Mathematics and Applied Mathematics\\ \small University of Pretoria\\\small Hatfield 0028, Pretoria, South Africa }
\date{\today}

\begin{document}
\maketitle

\begin{abstract}
We develop methods for determining key properties (simplicity and the dimension of radical) of flip subalgebras in Matsuo algebras. These are interesting classes of commutative non-associative algebras that were introduced within the broader paradigm of axial algebras.
\end{abstract}

\noindent {{\em Key words and phrases\/}: Axial algebra, Non-associative algebra, Matsuo algebra, radical, flip, flip subalgebra.
\hfill}\\

\section{Introduction}

The class of axial algebras was introduced by Hall, Rehren and Shpectorov in 2015 in \cite{HRS1,HRS2}. These are commutative non-associative algebras generated by special idempotents called axes. The adjoint action of an axis is semisimple with prescribed eigenvalues. Furthermore, products of eigenvectors with respect to the adjoint action obey a fusion law. (See the detailed definitions in the next section.) Axial algebras for graded fusion laws have large automorphism groups and hence they have direct relation to groups. The key example of this is the 196884-dimensional Griess algebra for the Monster sporadic simple group. In fact, the entire paradigm of axial algebras grew out of attempts to generalise this example. Outside of algebra, axial algebras arise in vertex operator algebras in mathematical physics and, more recently, in several other areas, such as non-linear solutions of PDEs (e.g. Hsiang algebras) in analysis \cite{T} and algebras of metric curvature tensors and flows in differential geometry \cite{DF}.

Matsuo algebras, corresponding to $3$-transposition groups, were initially introduced by Matsuo in 2003 in \cite{M} and then broadly generalised in \cite{HRS2}. They belong in the class of algebras of Jordan type $\eta$, i.e., their fusion law is similar to Peirce decomposition of Jordan algebras, where $\eta=\frac{1}{2}$. They constitute a very interesting class of
algebras whose properties largely remain unexplored. Furthermore, the recent double axis construction by Galt, Joshi, Mamontov, Shpectorov and Staroletov \cite{GJMSS} takes as input a Matsuo algebra and its involutory automorphism (a flip) and it outputs a flip subalgebra, which is an axial algebra of Monster type, i.e., its fusion law is similar to the one for the Griess algebra. There is currently a project underway to try to understand flip subalgebras and their properties.

One of the fundamental properties of an algebra is whether it is simple or whether it contains proper non-zero ideals. For Matsuo algebras, effective methods for deciding simplicity were developed by Hall and Shpectorov in \cite{HS}. In particular, it can be shown that connected Matsuo algebras are simple for almost all values of the parameter $\eta$. The exceptional values of $\eta$, called \emph{critical values}, are related by a simple formula to the eigenvalues of the collinearity matrix of the underlying Fischer space. The goal of the present paper is to further generalise these methods so that they apply to flip subalgebras.

For this purpose, we introduce, for a given flip subalgebra $A$, the concept of its ambient Matsuo algebra $M$, which can be understood as the smallest Matsuo algebra containing $A$. We show that if $A$ is non-simple, e.g. it has a non-zero radical $R(A)$, then the ambient algebra $M$ also has non-zero radical $R(M)$. More precisely, we have the following.

\begin{theorem} \label{main critical}
If $A$ is a flip subalgebra and $M$ its ambient Matsuo algebra then $R(A)=R(M)\cap A$.
\end{theorem}

That is, the critical values for $A$ and $M$ are essentially the same. Indeed, if $R(A)\neq 0$ then, clearly, also $R(M)\neq 0$. As for the converse, in principle, we can have a situation where $R(M)\neq 0$ while $R(A)=R(M)\cap A=0$, but this is quite rare.

In this way, we reduce the problem for flip subalgebras to the easier case of Matsuo algebras. After this, we focus on the dimension of $R(A)$, which we re-interpret, following \cite{KMS}, as the multiplicity of the eigenvalue $0$ for the Gram matrix of the Frobenius form on $A$ (a symmetric bilinear form that associates with the algebra product). We first handle the case where the $3$-transposition group $G$ for $M$ is irreducible. In this case, we develop a system of three linear equations, which allow us to find not just the multiplicity of $0$, but the multiplicities of all eigenvalues of the Gram matrix. Once this case is fully understood, we turn to the reducible case, which occurs when $G$ has a non-central normal $2$- or $3$-subgroup. This translates into a non-trivial invariant partition of the set of axes and eventually allows us to reduce the calculation of the eigenspace dimensions to the case of smaller $3$-transposition groups. This is the content of Corollaries \ref{main tau} and \ref{main theta}, which are somewhat technical and require additional $3$-transposition group notation. Hence we will not present them here.

Development of the theory is illustrated by examples of known flip subalgebras. These were taken from the upcoming papers on flip subalgebras coming from specific classes of Matsuo algebras (or rather the corresponding $3$-transposition groups). There are
currently several of such papers in progress, including by Joshi, Shpectorov and Shi \cite{JSS} on the symplectic and orthogonal (over $\F_2$) groups of $3$-transpositions, by Hoffman, Rodrigues and Shpectorov \cite{HRS3} on the unitary $3$-transposition groups, by Gorshkov and Staroletov \cite{GS} on the triality and Fischer groups, and finally, by Alsaaedi, Bovdi and Shpectorov \cite{ABS} on the first reducible cases, for the groups $2^{n-1}:S_n$ and $3^{n-1}:S_n$. These generalise the earlier work of Joshi on $S_n$ (see also, \cite{GJMSS}). All these projects are held back by the lack of methods of checking simplicity and computing the radical. Now that we have developed such methods in this paper, hopefully all these projects will be brought to completion.

The paper is organised as follows. In Section 2, we provide the necessary background on axial algebras and, in Section 3, we, similarly, introduce $3$-transposition groups and related notation. In Section 4, we define Matsuo algebras, which are the key objects we study. Section 5 discussed the concepts of radical and connectivity in axial algebras. Here we also review how these these can be determined in Matsuo algebras. In Section 6, we define flips and flip subalgebras and introduce the concept of the ambient Matsuo algebra, which is essential for the later considerations. In Section 7, we prove that the flip subalgebras have essentially the same critical values of $\eta$ as their ambient Matsuo algebras (Theorem \ref{main critical}). In Section 8, we develop linear equations which allow us to determine the dimension of the radical in the irreducible case. This approach is illustrated by two examples in Section 9. The reducible case is the focus of Section 10, while in Section 11 we illustrate the reduction methods on some examples. Throughout the paper we tacitly assume that the ground field is of characteristic $0$,
hence we also include, in Section 12, a brief discussion of what happens in positive characteristic.

\section{Background on axial algebras}

We start by reviewing the basics of axial algebras. For a more detailed discussion, see \cite{MS}.

\begin{definition}
A \emph{fusion law} is a set $\cF$ together with a map $\ast:\cF\times\cF\to 2^\cF$, where
$2^\cF$ denotes the set of all subsets of $\cF$.
\end{definition}

Normally, we think of $\cF$ as a small finite set and represent the map $\ast$ by a table similar to the multiplication
table of, say, a group. For example, Table \ref{fusion laws} shows the two fusion laws. that will appear later in
this paper.
\begin{figure}[h]
\setlength{\tabcolsep}{4pt}
\renewcommand{\arraystretch}{1.5}
\centering
	\begin{minipage}[t]{0.25\linewidth}
		\begin{tabular}{|c||c|c|c|}
		\hline
		$\star$ & $1$ &$0$&$\eta$\\
		\hline  \hline
		$1$&$1$& &$\eta$\\
		\hline
		$0$& &$0$&$\eta$\\
		\hline
		$\eta$&$\eta$&$\eta$&$1,0$\\
		\hline
		\end{tabular}
	\end{minipage}
	\begin{minipage}[t]{0.30\linewidth}
		\begin{tabular}{|c||c|c|c|c|}
		\hline
		$\star$ & $1$ &$0$&$\alpha$& $\beta$\\
		\hline \hline
		$1$&$1$& &$\alpha$& $\beta$\\
		\hline
		$0$& &$0$&$\alpha$& $\beta$\\
		\hline
		$\alpha$&$\alpha$&$\alpha$&$1,0$& $\beta$\\
		\hline
		$\beta$&$\beta$&$\beta$&$\beta$&$1,0, \alpha$\\
		\hline
		\end{tabular}
	\end{minipage}
\caption{The $\cJ(\eta)$, $\cM(\alpha,\beta)$ fusion laws.}
\label{fusion laws}
\end{figure}

Now suppose that $A$ is a commutative non-associative algebra over a field $\F$. For an element $a\in A$, we denote by
$\ad_a$ the \emph{adjoint map} of $a$, that is, the linear map $A\to A$ given by $u\mapsto au$. For $\lm\in\F$, let
$$A_\lm(a)=\{u\in A\mid au=\lm u\}$$
be the $\lm$-eigenspace of $\ad_a$. Note that $A_\lm(a)\neq 0$ if and only if $\lm$ is an eigenvalue of $\ad_a$. For
$\Lm\subseteq\F$, we write
$A_\Lm(a)=\oplus_{\lm\in\Lm}A_\lm(a)$.

\begin{definition}
An \emph{axis} in $A$ for the fusion law $\cF\subseteq\F$ is a non-zero idempotent $a\in A$ such that
\begin{enumerate}
\item[(a)] $A=A_\cF(a)$; and
\item[(b)] for $\lm,\mu\in\cF$, we have that $A_\lm(a)A_\mu(a)\subseteq A_{\lm\ast\mu}(a)$.
\end{enumerate}
\end{definition}

Condition (a), tells us that the linear map $\ad_a$ is semisimple, and furthermore, all eigenvalues of $\ad_a$ are contained
in the fusion law $\cF$. Condition (b) is the main part of the definition, and it means that the fusion law controls products
of eigenvectors of $\ad_a$.

From now on, we always assume that $cF\subseteq\F$. Furthermore, since $\ad_a(a)=aa=a$, we have that $1$ is an
eigenvalue of $\ad_a$, and so we must have $1\in\cF$.

\begin{definition}
An axis $a\in A$ is called \emph{primitive} if $A_1(a)=\la a\ra$ is $1$-dimensional.
\end{definition}

Now we are ready for the main concept.

\begin{definition}
Suppose $A$ is a commutative non-associative algebra over a field $\cF$ and $\cF\subseteq\F$ is a fusion law. Let
$X$ be a set of non-zero idempotents from $A$. We say that the pair $(A,X)$ is a (primitive) axial algebra for the fusion
law $\cF$ if
\begin{enumerate}
\item[(a)] every $a\in X$ is a (primitive) axis in $A$ for the fusion law $\cF$; and
\item[(b)] $X$ generates $A$, that is, $A=\lla X\rra$.
\end{enumerate}
\end{definition}

We note that the fusion law is used to distinguish classes of axial algebras. For example, primitive axial algebras for the
fusion law $\cJ(\eta)$ in Figure \ref{fusion laws} are called algebras of Jordan type $\eta$, whereas primitive axial algebras
for the fusion law $\cM(\al,\bt)$ in the same figure are called the algebras of Monster type $(\al,\bt)$. These are the two
classes of axial algebras that received the most attention in recent research.

Often an axial algebra admits a natural bilinear form.

\begin{definition}
A non-zero bilinear form $(\cdot,\cdot)$ on an axial algebra $(A,X)$ is called a \emph{Frobenius form} if it associates
with the algebra product, that is,
$$(uv,w)=(u,vw)$$
for all $u,v,w\in A$.
\end{definition}

An algebra admitting a Frobenius form is called \emph{metrisable}. For example, it was shown in \cite{HSS2} that every
algebra of Jordan type admits a unique Frobenius form such that $(a,a)=1$ for each generating $a\in X$. That is, all
algebras of Jordan type are automatically metrisable. All known examples of algebras of Monster type also admit a Frobenius
form, but there is currently no general result about metrisability of algebras of Monster type.

We note that all examples in this paper will be metrisable and, additionally, it will always be the case that $(a,a)\neq 0$ for
all $a\in X$. Then the form is typically scaled so that $(a,a)=1$ for all or at least some of the axes $a\in X$.

The structure theory for axial algebras was developed in \cite{KMS}.

\begin{definition}
The \emph{radical} of an axial algebra $(A,X)$ is the unique largest ideal not containing any axes from $X$.
\end{definition}

We note the following alternative description of the radical.

\begin{theorem} \label{radical}
If $A$ admits a Frobenius form such that $(a,a)\neq 0$ for each $a\in X$ then the radical of $A$ coincides with the radical
$$A^\perp=\{u\in A\mid (u,v)=0\mbox{ for all }v\in A\}$$
of the Frobenius form on $A$.
\end{theorem}

It is clear from the definition that the radical of $A$ contains all ideals of $A$ not containing any of the generating axes.
Can there be proper ideals containing generating axes? In general, this is possible, for example when $A$ is a direct sum of two
axial algebras. To control such ideals, we introduce the following concepts.

\begin{definition}
Suppose that $(A,X)$ is a (primitive) axial algebra admitting a Frobenius form such that $(a,a)\neq 0$ for all $a\in X$.
The \emph{projection graph} of $A$ is the graph on the set $X$, where distinct $a,b\in X$ are connected by an edge if and
only if $(a,b)\neq 0$.

We further say that $A$ is \emph{connected} if its projection graph is connected.
\end{definition}

The following result shows the consequences of connectivity.

\begin{theorem}
If $A$ admits a Frobenius form such that $(a,a)\neq 0$ for each $a\in X$ and $A$ is connected then $A$ has no proper
ideal containing generating axes. In particular, all proper ideals of $A$ are in the radical $A^\perp$.
\end{theorem}

We will now introduce a specific class of axial algebras, called Matsuo algebras, that will be the focus of this paper. The
definition of a Matsuo algebra involves a $3$-transposition group.

\section{$\mathbf{3}$-Transposition groups} \label{3-transposition}

Let us review the basic facts about $3$-transposition groups.

\begin{definition}
Suppose that $G$ is a group and $D\subset G$ is a normal subset of $G$. We say that
$(G,D)$ is a \emph{$3$-transposition group} if
\begin{enumerate}
\item[(a)] $D$ generates $G$; i.e., $G=\la D\ra$;
\item[(b)] all $d\in D$ have order $2$; and
\item[(c)] for all $d,e\in D$, the order of $de$ is at most $3$.
\end{enumerate}
\end{definition}

It is easy to see that $|de|=1$ if and only if $d=e$; $|de|=2$ if and only if $d$ and $e$
commute, but $d\neq e$; and $|de|=3$ if and only if $d$ and $e$ do not commute. In
the latter case, $H:=\la d,e\ra\cong S_3$, and so there is a unique additional element
$f\in H$ of order $2$. Since $d,e,f$ are conjugate in $H$, they are also conjugate in $G$ and
so $f\in D$. We have that $f=d^e=e^d$, $e=d^f=f^d$, etc.

This allows us to introduce a geometry $\Gm=\Gm(G,D)$ called the \emph{Fischer space}
related to the $3$-transposition group $(G,D)$. Its point set is $D$ and the lines are all triples
$\{c,d,e\}\subset D$, where $c$ and $d$ do not commute (i.e., $|cd|=3$) and $e=c^d=d^c$.
From the above, it is easy to see that this definition of line is symmetric with respect to the points
on it. In particular, any two collinear (non-commuting) points are contained in a unique line,
that is, the Fischer space is a partial linear space.

Clearly, two points on a line are conjugate in $G$, and so if $D$ is a union of several conjugacy
classes, $D=D_1\cup D_2\cup\ldots\cup D_k$, then every line is fully contained in some conjugacy
class $D_i$. This means that the Fischer space $\Gm$ is connected if and only if $D$ is
a single conjugacy class.  In the disconnected case, involutions from different classes
$D_i$ necessarily commute, which means that $G$ is the central product of the connected
$3$-transposition groups $(G_i,D_i)$, where $G_i=\la D_i\ra$ for $i=1,2,\ldots,k$. Because of
this reduction to the connected case, we will typically work with connected $3$-transposition
groups.

The following result was obtained by Fischer \cite{F}. Here we say that groups $G$ and $H$ have
the same \emph{central type} if $G/Z(G)\cong H/Z(H)$ as $3$-transposition groups.

\begin{theorem}
Suppose that $(G,D)$ is a non-abelian connected $3$-transposition group such that every soluble
normal subgroup of $G$ is contained in the centre $Z(G)$. Then $G$ is one of the following central
types:
\begin{enumerate}
\item[(a)] $S_n$, $n\geq 5$;
\item[(b)] $SO^\eps(2m,2)$, $m\geq 4$ if $\eps=+$ and $m\geq 3$ if $\eps=-$;
\item[(c)] $Sp(2m,2)$, $m\geq 3$;
\item[(d)] $SU(n,4)$, $n\geq 4$;
\item[(e)] $SO^\eps(n,3)$, $n\geq 6$;
\item[(f)] the triality groups $\Omega^+(8,2):S_3$, $\Omega^+(8,3):S_3$;
\item[(g)] sporadic Fischer groups, $Fi_{22}$, $Fi_{23}$, and $Fi_{24}$.
\end{enumerate}
\end{theorem}

Note that the only abelian connected group $(G,D)$ is of order $2$.

Cuypers and Hall \cite{CH} completed the classification of all connected $3$-transposition
groups. The statement of their classification is somewhat complicated, with many cases, and so we only explain
how an arbitrary connected finite $3$-transposition group can be reduced to one of the central
types in Fischer's theorem above.

For a connected $3$-transposition group $(G,D)$ and $d\in D$, define the sets
$A_d=\{e\in D\mid |de|=3\}$ and $D_d=\{e\in D\mid |de|=2\}$. We further introduce two
equivalence relations, $\tau$ and $\theta$, on $D$: we have that $d\tau e$ if and only if
$A_d=A_e$ and, similarly, $d\theta e$ if and only if $D_d=D_e$. This leads to two normal
subgroups: $\tau(G)=\la de\mid d\tau e, d,e\in D\ra$ and $\theta(G)=\la de\mid d\theta e,
d,e\in D\ra$. We note that only one of these two normal subgroups can be non-trivial.
Furthermore, $\tau(G)$ is a $2$-subgroup and $\theta(G)$ is a $3$-subgroup. If both
subgroups are trivial, i.e, $\tau(G)=1=\theta(G)$, then the $3$-transposition group
$(G,D)$ is called \emph{irreducible}. The irreducible groups are the ones listed in Fischer's
Theorem above plus the group of order $2$.

We are going to make use of the following notation. Suppose that $\tau(G)\neq 1$ then for
every $d\in D$, we have that $|D\cap d\tau(G)|=2^h$ for some $h\geq 1$, independent
of $d$. We also have the factor $3$-transposition group $(\bar G,\bar D)$, where
$\bar G=G/\tau(G)$ and $\bar D$ is the image of $D$ in $\bar G$. Then we have that
$|D|=2^h|\bar D|$. We will say that $(G,D)$ is a $2^h$-fold cover of $(\bar G,\bar D)$
and write $G=2^{\bullet h}\bar G$. (We note that $2^h$ is not the size
of $\tau(G)$, which can be significantly larger than $2^h$.) In some cases, we can write
$G=4^{\bullet h}\bar G$, and then, of course, $|D\cap d\tau(G)|=4^h$.

Similarly, when $\theta(G)\neq 1$, we have that $|D\cap d\theta(G)|=3^h$ and again
we will speak of $G$ as the $3^h$-fold cover of $\bar G=G/\theta(G)$ and write
$G=3^{\bullet h}\bar G$. Again, in some cases we can write $9^{\bullet h}\bar G$,
where $h$ is the half of its earlier counterpart.

\begin{example} \label{tau}
Get $\hat G=2^n:S_n$, the semidirect product of $G=S_n$ and its permutation module $V$
over the field with two elements. We will view the normal subgroup $V=2^n$ as the direct product
of $n$ copies of $C_2=\{1,-1\}$, the multiplicative version of the group of order $2$,
and correspondingly, we will write $(\dl_1,\dl_2,\ldots,\dl_n:\sg)$ for a general element
of $\hat G$. Here $\sg\in G$ is a permutation and $\dl_1,\dl_2,\ldots,\dl_n\in C_2=\{1,-1\}$.
For example, the element $\sg\in G$ can also be written as $(1,1,\ldots,1:\sg)\in\hat G$.

The class $D$ of $3$-transpositions in $G=S_n$ is the class of transpositions ($2$-cycles)
and it extends to the class $\hat D$ of $3$-transpositions in $\hat G$, which includes all
elements $(i,j)=(1,1,\ldots,1:(i,j))$ as well as all elements $v_{i,j}(i,j)=(1,\ldots, -1,\ldots, -1,\ldots,1:(i,j))$,
where $v_{i,j}=(1,\ldots,-1,\ldots,-1,\ldots,1:())$ has $-1$s in positions $i$ and $j$.
So we see that $|\hat D|=2|D|$ and so we can write the $3$-transposition subgroup $G^\circ$
generated by $\hat D$ inside $\hat G$ (this subgroup $G^\circ$ has index $2$ in $\hat G$) as
$2^{\bullet 1}S_n$.

Indeed, in this case $d\tau d'$ for $d,d'\in\hat D$, $d\neq d'$, if and only if $\{d,d'\}=
\{(i,j),v_{i,j}(i,j)\}$ for some $1\leq i<j\leq n$. Hence $\tau(G^\circ)=\la v_{i,j}\mid 1\leq i
<j\leq n\}$ and this is an index $2$ subgroup $V^\circ$ of $V$. We then have that
$G^\circ=V^\circ G\cong 2^{n-1}:S_n$ and for $d\in\hat D$, we have that
$\hat D\cap d\tau(G^\circ)=\{d,d'\}=\{(i,j),v_{i,j}(i,j)\}$. So we do indeed have $h=1$ in
this example.
\end{example}

Similarly, if we consider $\hat G=2^{2m}:Sp(2m,2)$, the semidirect product of $G=Sp(2m,2)$
with its natural module, then this will also have $h=1$ and so we can write $2^{\bullet 1}Sp(2m,2)$
for this $3$-transposition group.

\begin{example}
The second example shows the case where $\theta$ is a non-trivial relation. Let $\hat G=3^n:S_n$
be the semidirect product of $G=S_n$ with its permutation module $V$ over the field with three
elements. We will use the notation similar to Example $\ref{tau}$, namely we will write a general
element of $\hat G$ as $(\dl_1,\dl_2,\ldots,\dl_n:\sg)$, where $\sg\in S_n$ and $\dl_1,\dl_2,\ldots,
\dl_n\in C_3=\{1,\xi,\xi^{-1}\}$, where $\xi=e^{2\pi i/3}$, the multiplicative version of the
group of order $3$. Let $v_i:=(1,\ldots,\xi,\ldots,1:())$, with $\xi$ in position $i$. Let also
$v_{i,j}:=v_iv_j^{-1}$. (So that $v_{i,j}=v_{j,i}^{-1}$.). Then the $3$-transposition class
$\hat D$ extending the class $D$ from $G=S_n$, consists of all elements of the form $(i,j)$,
$v_{i,j}(i,j)$, and $v_{j,i}(i,j)$, for $1\leq i<j\leq n$. Furthermore, such triples from $\hat D$
(for fixed $i$ and $j$) are the equivalence classes of the relation $\theta$ on $\hat D$.
From this description, it is clear that $\theta(\hat G)=\la v_{i,j}\mid 1\leq i<j\leq n\ra$ has index
$3$ in $V$ and, in fact, $\hat D$ generates in $\hat G$ a subgroup of index $3$ isomorphic to
$\theta(\hat G)S_n$.

Since the equivalence classes of $\theta$ on $\hat D$ have size $3$, we have that this $3$-transposition
group can be written as $3^{\bullet 1}S_n$.
\end{example}

We also get similar examples from other cases of irreducible $3$-transposition groups. We will see some
of them later.

\section{Matsuo algebras} \label{Matsuo}

At the centre of this paper lies the class of algebras known as Matsuo algebras. Consider a
$3$-transposition group $(G,D)$. Recall that the Fischer space of $(G,D)$ has $D$ as the set of points
and it has as lines the set of all triples of points $\{c,d,e\}\subseteq D$, such that $|cd|=3$ and
$e=c^d=d^c$.

Select a field $\F$ of characteristic not $2$ and $\eta\in\F$, $\eta\notin\{1,0\}$.
Define the Matsuo algebra $M=M_\eta(G,D)$ as follows: We take $D$ to be the basis
of $M$ and, for $c,d\in D$, we define
$$cd=\left\{\begin{array}{ll}
c,&\mbox{if }c=d;\\
0,&\mbox{if }|cd|=2;\\
\frac{\eta}{2}(c+d-c^d),&\mbox{if }|cd|=3.
\end{array}
\right.
$$

The Matsuo algebra $M$ is an algebra of Jordan type $\eta$, namely, as shown in \cite{HRS2},
the elements of the basis $D$ of $M$ are primitive axes of Jordan type $\eta$ and they clearly
generate $M$.

Note that we did not require above that $(G,D)$ be connected. If $(G,D)$ is disconnected, i.e.,
if $D=D_1\cup D_2\cup\ldots\cup D_k$ is a union of $k>1$ conjugacy classes, then $c$ and $d$
from different classes necessarily commute in $G$, that is, they multiply to zero in $M$.
This immediately implies that $M=\oplus_{i=1}^k M_i$, where $M_i=M_\eta(G_i,D_i)$ is
a similar Matsuo algebra of the smaller connected $3$-transposition group $(G_i,D_i)$, with
$G_i=\la D_i\ra$, $i=1,2,\ldots,k$. Thus, the disconnected case simply leads to direct sums
and this allows us in most cases to restrict the problem to the connected case.

We note that the Matsuo algebra $M$ admits a Frobenius form, and its values on the basis $D$
are as follows:
$$(c,d)=\left\{\begin{array}{ll}
1,&\mbox{if }c=d;\\
0,&\mbox{if }|cd|=2;\\
\frac{\eta}{2},&\mbox{if }|cd|=3.
\end{array}
\right.$$

\section{Radical and connectivity} \label{radical and connectivity}

Our focus is on the properties of algebras. So let us first discuss connectivity
and radicals in Matsuo algebras. We start with connectivity.

\begin{proposition}
The projection graph of $M=M_\eta(G,D)$ coincides with the collinearity graph of the Fischer
space of $(G,D)$. In particular, $M$ is connected if and only if the Fischer space is connected
if and only if $D$ is a single conjugacy class of $G$.
\end{proposition}

\pf The first claim  follows directly from the definition of the projection graph, since for
distinct $c,d\in D$, we have that $(c,d)\neq 0$ if and only if $|cd|=3$, i.e., they are collinear
in the Fischer space.\qed

\medskip
We saw in the preceding section that in the disconnected case $M$ decomposes as the direct
sum of the smaller connected Matsuo algebras corresponding to the components of the
projection graph. For the remainder of this section, we assume that $M$ is connected. This
means that $M$ is simple if and only if it has zero radical and this depends, as we will soon
see, on the value of $\eta$.

\begin{definition} \label{critical}
Suppose that $(G,D)$ is a connected $3$-transposition group and $M=M_\eta(G,D)$ its
Matsuo algebra. We call the value $\eta$ \emph{critical} if $M$ is not simple, i.e., if $M$ has
a non-zero radical.
\end{definition}

Now recall from Theorem \ref{radical} that the radical of $M$ coincides with the radical
$M^\perp$ of the Frobenius form on $M$ and so it is non-zero if and only if the nullspace
of the Gram matrix $L$ of the Frobenius form on $M$ is non-zero. That is, if and only if $0$
is an eigenvalue of $L$. Naturally, we write $L$ with respect to the basis $D$.

Let $T=(t_{cd})_{c,d\in D}$ be the collinearity matrix of the Fischer space, that is,
$t_{cd}=1$ if $c$ and $d$ are collinear and distinct (i.e., $|cd|=3$); and $t_{cd}=0$ otherwise.
This is the same as the adjacency matrix of the projection graph of $M$. Looking at the values
of the Frobenius form on the basis $D$, shown in the preceding section, we have the following
expression for the Gram matrix $L$:
$$L=I+\frac{\eta}{2}T,$$
where $I$ is the identity matrix of size $|D|$. Clearly, this means that the eigenvalues of $L$
are of the form $1+\frac{\eta}{2}\rho$, where $\rho$ runs through the eigenvalues of $T$.

Note that the eigenvalue formula has already appeared in \cite{HS}, where the projection graph
above is called the diagram of the group $(G,D)$. In particular, in \cite{HS}, one can find
the complete tables of the eigenvalues $\rho$ for all connected $3$-transposition groups $(G,D)$.

The tables from \cite{HS} allow us to determine the critical values of each connected Matsuo
algebra $M=M_\eta(G,D)$. Namely, $\eta$ is critical if and only if $0$ is an eigenvalue of $L$,
which means that we must have $0=1+\frac{\eta}{2}\rho$ for some eigenvalue $\rho$ of $T$.
Therefore, we must have $\eta=-\frac{2}{\rho}$. Consequently, as $T$ has only finitely many
eigenvalues, $M=M_\eta(G,D)$ has a finite number of critical values. In other words, connected
Matsuo algebras are generically simple.

\begin{example}
Let $G=S_n$ and $D=(1,2)^{S_n}$ be the class of transpositions. According to \cite{HS}, the
collinearity matrix $T$ has in this case eigenvalues $(2(n-2)^1,(n-4)^{n-1},
-2^{\frac{n(n-3)}{2}})$. (Here the exponent indicates the multiplicity of the eigenvalue.)

Based on this, we obtain that the critical values for $M=M_\eta(S_n,(1,2)^{S_n})$ are
$\eta=-\frac{2}{2(n-2)}=-\frac{1}{n-2}$, $\eta=-\frac{2}{n-4}$, and $\eta=-\frac{2}{-2}=1$.
However, the latter value should be dropped because $\eta\notin\{1,0\}$.
\end{example}

We will see more examples later in the paper.

Let us also discuss the dimension of the radical of $M$ when it is non-trivial. Clearly, when
$\eta=-\frac{2}{\rho}$ is critical then the dimension of the radical of $M$ is the dimension
of the $0$-eigenspace of $L$ and the latter coincides with the $\rho$-eigenspace of $T$. Thus,
the dimension of the radical is simply the multiplicity of the eigenvalue $\rho$. For instance,
in the above example, the dimension of the radical of $M$ when $\eta=-\frac{1}{n-2}$ is
the multiplicity of the eigenvalue $\rho=2(n-2)$ of $T$, which is $1$. Similarly, if $\eta=
-\frac{2}{n-4}$ then the dimension of the radical of $M$ is the multiplicity of the eigenvalue
$\rho=n-4$ of $T$, which is $n-1$. This calculation applies for fields $\F$ of characteristic
zero and also in all positive characteristics for which the relevant (integer) eigenvalues of
$T$ (as in the tables from \cite{HS}) are distinct from each other.

To conclude, the tables in \cite{HS} allow us to determine the critical values of any Matsuo
algebra $M$ and hence check its simplicity. They also allow us to find the dimension of the radical
of $M$ when it is non-zero, provided that the characteristic of the field $\F$ is not bad for
the corresponding $3$-transposition group $(G,D)$.

\section{Flips and flip subalgebras}

Now that we understand the case of Matsuo algebras, let us introduce a second important class
of algebras, which were defined in \cite{GJMSS} as subalgebras of Matsuo algebras.

Let $M=M_\eta(G,D)$ be a connected Matsuo algebra. An automorphism of the
$3$-transposition group $(G,D)$ is an automorphism $\sg$ of the group $G$ preserving the
class $D$, that is, $D^\sg=D$. Clearly, such an automorphism also acts on the Matsuo algebra
$M$ by permuting its basis $D$. When $|\sg|=2$, we say that $\sg$ is a \emph{flip} of $M$.

Let $\Sg=\la\sg\ra\le\Aut(M)$ be the corresponding subgroup of order $2$ and consider its orbits
on the basis $D$. We distinguish three types of such orbits: (1) orbits $\{c\}$ of length one,
for all $c\in D$ satisfying $c^\sg=c$; (2) orbits $\{c,d\}$ of length two (i.e., $c^\sg=d$ and
$d^\sg=c$) where $cd=0$, i.e., $c$ and $d$ are \emph{orthogonal}; and (3) orbits $\{c,d\}$ of length two, where $cd\neq 0$, i.e., $c$ and $d$ are non-orthogonal. We denote the set of all orbits of length one $O_s$, the set of all orthogonal orbits of length two $O_d$, and the set of non-orthogonal orbits of length two $O_e$.

For each orbit $O$, we define the sum vector $v_O=\sum_{c\in O}c\in M$. Clearly, if $O=\{c\}$
has length one then $v_O=c$ is an axis in $M$, and we will call such orbit sums \emph{singles}.
If $O=\{c,d\}$ is of second kind then $v_O=c+d$ is also an idempotent since $(c+d)^2=
c^2+cd+dc+d^2=c+d$, as $cd=dc=0$. We call such idempotents \emph{doubles}. Finally, if
$O=\{c,d\}$ is of the third kind then $v_O=c+d$ is not an idempotent at all. Indeed,
$(c+d)^2=c^2+cd+dc+d^2=c+d+2cd\neq c+d$. (Recall that the ground field cannot be of
characteristic $2$.) We call such vectors $v_O=c+d$ \emph{extras}.
Thus, $|O_s|$, $|O_d|$, and $|O_e|$ are the numbers of singles, doubles, and extras, respectively.

Note that all singles, doubles, and extras form a basis in the fixed subalgebra
$$M_\sg=\{u\in M\mid u^\sg=u\}.$$

The following key observation was made in \cite{GJMSS} (see Propositions 4.7 and 7.4 there).

\begin{theorem}
If $\eta\neq\frac{1}{2}$ then the doubles are axes of Monster type $(2\eta,\eta)$. They are
imprimitive in $M$, but primitive in $M_\sg$.
\end{theorem}

This naturally leads to the following concept.

\begin{definition}
The \emph{flip subalgebra} of $M=M_\eta(G,D)$ corresponding to the flip $\sg\in\Aut(M)$
is the subalgebra $A(\sg)$ of $M_\sg$ generated by all singles and doubles.
\end{definition}

As follows from the above discussion, $A(\sg)$ is an algebra of Monster type $(2\eta,\eta)$.
It can be viewed as twisted version of the Matsuo algebra $M$. Note that a $3$-transposition
group $(G,D)$ has many conjugacy classes of automorphisms of order two, and consequently,
the corresponding Matsuo algebra $M=M_\eta(G,D)$ has many classes of flip subalgebras.
Also note that, while the extras play no role in the definition of $A(\sg)$, they may in some
cases end up in $A(\sg)$, if they can be obtained by adding and multiplying
singles and doubles.

For the remainder of this paper, we study how the methods giving us the properties of Matsuo
algebras can be generalised to cover also the class of flip subalgebras.  However, the first order
of business is to notice that $M$ may not be the best algebra to construct $A(\sg)$ from.

\begin{definition}
Suppose $M=M_\eta(G,D)$ is a Matsuo algebra and let $A=A(\sg)$ be the flip subalgebra
corresponding to the flip $\sg\in\Aut(M)$. Let $C$ be the set of axes $c\in D$ involved in the
singles and doubles with respect to $\sg$. That is,
$$C=\{c\in D\mid c^\sg=c\mbox{ or }c^\sg c=0\}.$$
Let $(G',D')$ be the smallest $3$-transposition subgroup of $(G,D)$ containing $C$.
Then $(G',D')$ is called the \emph{ambient $3$-transposition group} of $A$ and
$M'=M_\eta(G_0,D_0)$ is called the \emph{ambient Matsuo algebra} of $A$.
\end{definition}

This definition requires a brief justification. Note that $(G',D')$ is minimal subject to:
(a) $C\subseteq D'\subseteq D$; (b) $G'=\la D'\ra$; and (c) $D'$ is a normal
subset of $G'$. This 3-transposition subgroup $(G',D')$ can be constructed as follows:
$G'=\la C\ra$ and $D'=C^{G'}$. Furthermore, while we cannot claim in general
that $D'=C$, this in fact will be the case in the examples we will see later.

Let us now see how this concept fits into the paradigm of flip subalgebras.

\begin{theorem}
Suppose $M_\eta(G,D)$ is the Matsuo algebra for $(G,D)$. Suppose further that $\sg$ is
a flip of $M$ and $A=A(\sg)$ is the corresponding flip subalgebra. Then
\begin{enumerate}
\item[(a)] the ambient subalgebra $M'=M_\eta(G',D')$ of $A$ is invariant under $\sg$
(i.e., $\sg$ is also a flip of $M'$); and
\item[(b)] the flip subalgebra of $M'$ coming from the flip $\sg$ coincides with $A$.
\end{enumerate}
\end{theorem}

\pf Clearly $\sg$ leaves $C$ invariant, and so $\sg$ normalises $G'=\la C\ra$ and
$D'=C^{G'}$. Thus, $\sg$ indeed induces a flip of $(G',D')$. Furthermore,
since $C\subseteq D'$, the singles and the doubles arising from the action of $\Sg=\la\sg\ra$
on $D'$ are the same as the ones arising from the action of $\Sg$ on $D$. Hence
the flip subalgebra of $M'$ coming from $\sg$ coincides with $A$.\qed

\medskip
Hence we can always reduce to the case where $M$ is the ambient Matsuo subalgebra of
$A$. We will see in the examples below that, under this additional assumption, $A(\sg)$
often coincides with the entire fixed subalgebra $M_\sg$.

\section{Critical values}

We call a value of $\eta\in\F$ \emph{critical} for a flip subalgebra $A=A(\sg)$ of
$M=M_\eta(G,D)$ if $A$ has a non-zero radical. Clearly, this is a straightforward
generalisation of the corresponding concept for connected Matsuo algebras introduced in Definition
\ref{critical}. (We will also use this term in a similar way for other subalgebras of Matsuo algebras, e.g. the fixed subalgebras $M_\sg$.) In this section we detail how the critical values of $A$ can be
determined depending on $(G,D)$ and $\sg$. Recall that, as we have already mentioned,
typically the flip subalgebra $A$ coincides with the entire fixed subalgebra $M_\sg$,
provided that $M$ is chosen to be the ambient Matsuo algebra of $A$. So we just need
to see for which values of $\eta$ the fixed subalgebra $M_\sg$ has a non-zero radical.

In addition to $M_\sg$, let us consider the commutator space
$$[M,\sg]=\{u-u^\sg\mid u\in M\}.$$
While $M_\sg$ has a basis consisting of all singles, doubles and extras, i.e., the orbit sums
for the group $\Sg=\la\sg\ra$, the commutator subspace $[M,\sg]$ has a spanning set consisting
of the commutators,  $C=\{[a,\sg]=a^\sg-a\mid a\in D\}$. Note that if $a$ is a single then
$a^\sg-a=a-a=0$, so we can remove all such zero vectors from $C$. Furthermore, if
$\{a,b\}=\{a,a^\sg\}\in O_d$ or $O_e$ is an orbit of length two then the commutators $a^\sg-a=b-a$ and $b^\sg-b=a-b$
add up to zero, and so we can leave in $C$ just one commutator for each such orbit. This modified set
$C^\circ$, without zero vectors and multiples, is linearly independent, because now the vectors from $C^\circ$
have pairwise disjoint support, and so it is a basis of $[M,\sg]$. We will call the vectors $a^\sg-a$ for
orthogonal orbits $\{a,a^\sg\}\in O_d$ \emph{skew doubles}, and for non-orthogonal orbits
$\{a,a^\sg\}\in O_e$ \emph{skew extras}.

\begin{proposition} \label{orthogonal}
The commutator space $[M,\sg]$ is orthogonal to $M_\sg$ with respect to the Frobenius form
on $M$.
\end{proposition}

\pf Let $u\in M_\sg$, that is, $u^\sg=u$, and let $w=[v,\sg]=v-v^\sg\in[M,\sg]$. Then $(u,w)=
(u,v-v^\sg)=(u,v)-(u,v^\sg)=(u,v)-(u^\sg,v^\sg)=0$.\qed

\medskip
For a subspace $U$ of $M$, we define the \emph{radical} of $U$ as
\begin{align*}
R(U)&=U\cap U^\perp=U\cap\{v\in A\mid (u,v)=0\mbox{ for all }u\in U\}\\
&=\{v\in U\mid (u,v)=0\mbox{ for all }u\in U\}.
\end{align*}

\begin{proposition}
We have that $M=M_\sg\oplus [M,\sg]$ is an orthogonal direct sum and, consequently, $M^\perp=R(M)
=R(M_\sg)\oplus R([M,\sg])$.
\end{proposition}

\pf Clearly, $\dim M=|D|=|O_s|+2|O_d|+2|O_e|=(|O_s|+|O_d|+|O_e|)+(|O_d|+|O_e|)=\dim M_\sg+
\dim [M,\sg]$. Also, by observation, $M_\sg+[M,\sg]=M$, which means that $M=M_\sg\oplus
[M,\sg]$. By Proposition \ref{orthogonal}, $M_\sg\perp [M,\sg]$, and so this direct sum is
orthogonal.

Clearly, $R(M_\sg)\oplus R([M,\sg])$ is orthogonal to $M_\sg\oplus [M,\sg]=M$, that is,
$R(M_\sg)\oplus R([M,\sg])\subseteq M^\perp$. On the other hand, suppose $v\in M^\perp$ and
$v=u+w$, where $u\in M_\sg$ and $w\in [M,\sg]$. Then $w$ and $v$ are orthogonal to $M_\sg$ and
so $u=v-w$ is also orthogonal to $M_\sg$, that is, $u\in R(M_\sg)$. Symmetrically, $w\in
R([M,\sg])$.\qed

\medskip
As a consequence, we get the following.

\begin{corollary} \label{radical dimension}
$\dim R(M_\sg)=\dim R(M)-\dim R([M,\sg])$.
\end{corollary}

In particular, if $R(M_\sg)\neq 0$ then also $R(M)\neq 0$. Therefore, we can claim the
following.

\begin{theorem}
The critical values of $M_\sg$ are also critical for $M$.
\end{theorem}

As we already mentioned, the flip subalgebra $A$ usually coincides with $M_\sg$ provided that $M$ is
the ambient Matsuo algebra for $A$. Hence, the above theorem tells us what the critical values of the
majority of flip subalgebras are.

Clearly, the multiplicity of each critical value for $M_\sg$ is going to be smaller than its multiplicity in $M$. We will view the situation where a critical value of $M$ is not critical for $M_\sg$ as simply the instance where the multiplicity becomes zero in $M_\sg$.

\section{Multiplicity}

In this section we look at the dimension of the radical of flip subalgebra $A=A(\sg)$ and we will again assume
that $A$ coincides with $M_\sg$. Recall that the radical $M_\sg^\perp$ of $M_\sg$ coincides with the radical of the Frobenius form on $M_\sg$
and hence the dimension of the radical is the same as the multiplicity of the eigenvalue zero for
the Gram matrix of the Frobenius form restricted to $M_\sg$.

Recall that $M$ is defined in terms of a connected $3$-transposition group $(G,D)$ and, furthermore, recall the notation $A_d=\{e\in D\mid |de|=3\}$ and $D_d=\{e\in D\mid |de|=2\}$. Let $t:M\to M$ be the linear
map taking every $d\in D$ to $t(d)=\sum_{c\in A_d}c$, the sum of all points $c\neq d$ of the Fischer
space that are collinear with $d$. We call $t$ the \emph{collinearity map}. Then the collinearity
matrix $T$ of the Fischer space $D$ is exactly the matrix of $t$ with respect to the basis $D$ of $M$. This means that we
can study the eigenvalues of $t$ instead of $T$.

\begin{proposition} \label{commute}
The maps $t$ and $\sg$ commute.
\end{proposition}

\pf It suffices to check that $t$ and $\sg$ commute when applied to the basis $D$. In
view of the above, for $d\in D$, we have that $t(d)^\sg=(\sum_{c\in A_d}c)^\sg=\sum_{c\in
A_d}d^\sg=\sum_{e\in A_d^\sg}e=\sum_{e\in A_{d^\sg}}e=t(d^\sg)$.\qed

\begin{corollary} \label{invariant}
The collinearity map $t$ leaves both $M_\sg$ and $[M,\sg]$ invariant.
\end{corollary}

\pf Indeed, as $t$ and $\sg$ commute, if $u\in M_\sg$ then $u^\sg=u$ and hence
$t(u)^\sg=t(u^\sg)=t(u)$. That is, $t(u)\in M_\sg$, proving that $M_\sg$ is invariant under
$t$.

Similarly, if $w\in [M,\sg]$ then $w=u^\sg-u$ for some $u\in M$. Therefore,
$t(w)=t(u^\sg-u)=t(u^\sg)-t(u)=t(u)^\sg-t(u)=[t(u),\sg]\in [M,\sg]$. Thus, $[M,\sg]$ is also
invariant under $t$.\qed

\medskip
It follows from Corollary \ref{radical dimension} that, in order to find the dimension of the radical of
$M_\sg$, it suffices to find the dimension of the radical of $[M,\sg]$, since
$\dim M^\perp=\dim R(M_\sg)+\dim R([M,\sg])$. The
advantage of dealing with $[M,\sg]$ is that its basis, consisting of skew doubles and skew extras,
is more uniform, as it only involves orbits of $\Sg$ of length two.

Note that, by Corollary \ref{invariant}, the map $t$ acts on $W$ and, clearly, the eigenvalues of
$t$ on $W$ must be within the set of eigenvalues of $t$ on $M$. This set is small and so, in order
to find the multiplicities of the eigenvalues of $t$ in $W$, we just need to derive a few linear equations
involving them.

Let $\Lm$ be the set of eigenvalues of $t$ on $M$, and let $m_\lm$, for $\lm\in\Lm$, be the
multiplicity of $\lm$ within $W$. (We note that $m_\lm$ may be zero; that is, not every
eigenvalue of $t$ arises on $W$.) The first equation simply expresses the fact that the sum
of all multiplicities is equal to $\dim W$.

\begin{proposition}
We have that
\begin{equation} \label{eq 1}
\sum_{\lm\in\Lm}m_\lm=\dim W.
\end{equation}
\end{proposition}

\pf Since the collinearity relation is symmetric, the matrix $T$ of $t$ with respect to the
basis $D$ of $M$ is symmetric. This implies that $T$ is diagonalisable, that is, its minimal
polynomial $p$ has no multiple roots.

This, in turn, means that the minimal polynomial of $t$ on $W$, being a divisor of $p$, also
has no multiple roots. Thus, $t$ must be diagonalisable. This means that $W$ is a direct sum
of eigenspaces of $t$, and so equation (\ref{eq 1}) must hold, since $m_\lm$ is the dimension
of the $\lm$-eigenspace of $t$.\qed

\medskip
The second equation comes from the consideration of the trace of $t$ on $W.$ Clearly, we have
that the trace equals $\sum_{\lm\in\Lm}m_\lm\lm$, since $t$ is diagonalisable on $W$. On
the other hand, the trace of $t$ does not depend on the choice of basis, and so we can find it
with respect to the basis consisting of the skew doubles and skew extras.

\begin{proposition}
The trace of $t$ on $W$ is equal to minus the number of skew extras in the basis of $W$.
\end{proposition}

\pf Let $w_1,w_2,\ldots,w_m$ be a basis of $W$ consisting of skew doubles and skew extras.
Take $w=w_i=d-d^\sg$ for some $d\in D$. Then $t(w)=t(d)-t(d^\sg)=t(d)-t(d)^\sg=\sum_{c\in
A_d}c-\sum_{c\in A_d}c^\sg=\sum_{c\in A_d}(c-c^\sg)$. (Some of the summands here can
be zero; this happens when $c=c^\sg$ is a single.) If $w$ is a skew double then $dd^\sg=0$
in $M$, and so $d^\sg\notin A_d$ (clearly, $d$ is never in $A_d$). This means that the $i$th
position in the $i$th column of the matrix of $t$ with the above basis is $0$, as $d-d^\sg$
does not appear in the linear decomposition of $t(d-d^\sg)$. So this basis vector does not
contribute to the trace.

On the other hand, if $w=w_i=d-d^\sg$ is a skew extra then $dd^\sg\neq 0$ in $M$; that is,
$d^\sg\in A_d$. Looking at the expression for $t(w)$, we find there $d^\sg-(d^\sg)^\sg
=-(d-d^\sg)$, which means that each such basis vector $w=w_i$ contributes $-1$ to the trace.
Now the claim follows.\qed

\begin{corollary}
We have that
\begin{equation} \label{eq 2}
\sum_{\lm\in\Lm}m_\lm\lm=-e,
\end{equation}
where $e$ is the number of non-orthogonal orbits of $\Sg=\la\sg\ra$ on $D$.
\end{corollary}

We note that one of the eigenvalues $\lm\in\Lm$ is $k=|A_d|$, $d\in D$, the valency of the
collinearity graph on $D$.

\begin{proposition}
We have that
\begin{equation} \label{eq 3}
m_k=0.
\end{equation}
\end{proposition}

\pf Indeed, the multiplicity of $k$ in $M$ is $1$ (since we only consider connected cases),
and the corresponding eigenvector can be chosen as $\sum_{d\in D}d$. Obviously, this
eigenvector is in $M_\sg$, which means that the multiplicity $m_k$ of $k$ in $W$ is $0$.\qed

\medskip
The three equations we derived are sufficient to find the multiplicities $m_\lm$ in all cases
where $|\Lm|=3$, which is always the case for irreducible $3$-transposition groups. We will
discuss the reducible case later, after we give examples of computing the multiplicities in
the irreducible case.

\section{Two examples}

As an example, let us consider the case of $G=S_n$. Here $D$ is the class of transpositions
($2$-cycles). The eigenvalues of the adjacency matrix $T$ of the collinearity graph on $D$ and
the multiplicities of these eigenvalues can be found, for example, in \cite{HS}. They are
$k=2(n-2)$ (with multiplicity $1$), $r=n-4$ (with multiplicity $n-1$), and $s=-2$ (with
multiplicity $\frac{n(n-3)}{2}$). Correspondingly, the critical values for $M=M_\eta(G,D)$ are
$\eta=-\frac{2}{k}=-\frac{2}{2(n-2)}=-\frac{1}{n-2}$, for which the radical of $M$ is
$1$-dimensional, $\eta=-\frac{2}{r}=-\frac{2}{n-4}$, for which the radical is of dimension
$n-1$ and $\eta=-\frac{2}{s}=-\frac{2}{-2}=1$, for which the radical would have been of
dimension $\frac{n(n-3)}{2}$. However, we need to exclude this value of $\eta$, because
$\eta\notin\{0,1\}$ (see Section \ref{Matsuo}).

Turning now to the flip subalgebras of $M$, they have been classified by Joshi in \cite{J}
(see also \cite{GJMSS}). Namely, for every $m=1,2,\ldots,\lfloor\frac{n}{2}\rfloor$, we have,
up to conjugation, a single flip $\sg_m$, and the corresponding flip subalgebra is isomorphic
to $Q_m(\eta)\oplus M_\eta(S_{n-2m})$, where $Q_m(\eta)$ is the flip subalgebra arising
from the same $\sg$ when $n=2m$. This subalgebra has dimension $m^2$ and its standard
basis consists of $m$ singles, $m(m-1)$ doubles, and no extras.

We can use our equations to find the critical values and radical dimensions for $Q_m$. The
ambient algebra for $Q_m$ is $M$ if we again assume that $n=2m$. Thus, as for $M$, the critical
values are $\eta=-\frac{1}{n-2}=-\frac{1}{2m-2}$, and $-\frac{2}{n-4}=-\frac{2}{2m-4}=
-\frac{1}{m-2}$. To find the dimension of the radical of $A=Q_m(\eta)$, we first find the
multiplicities $m_\lm$ of the corresponding eigenvalues $2(n-2)$, $n-4$, and $-2$ in
$W=[M,\sg_m]$. (Here we include $\lm=-2$, which corresponds to the prohibited critical
value $\eta=1$, because we need to account for all eigenvalues of $T$.)
First of all, since $k=2(n-2)$, we have that $m_{2(n-2)}=0$ by (\ref{eq 3}). Also,
(\ref{eq 1}) and (\ref{eq 2}) become
\begin{equation} \label{one}
m_{2m-4}+m_{-2}=\dim W=\dim M-m^2=m(2m-1)-m^2=m^2-m,
\end{equation}
and
\begin{equation} \label{two}
(2m-4)m_{2m-4}-2m_{-2}=0,
\end{equation}
since the number $e$ of non-orthogonal orbits of length $2$ is zero, as we have no extras when $n=2m$.

Expressing $m_{-2}$ from equation (\ref{two}), we get $m_{-2}=(m-2)m_{2m-4}$. Substituting this into
(\ref{one}) yields $(m-1)m_{2m-4}=m^2-m$, which implies that $m_{2m-4}=m$. This in turn gives
$m_{-2}=m(m-2)$.

Now we can find the dimensions of the eigenspaces in $A=M_\sg=Q_m(\eta)$.

\begin{theorem} \label{orthogonal example}
Suppose that $A=M_\sg=Q_m(\eta)$. Then the radical of this algebra is non-zero if and only if $\eta$ is
in the table below, where in the bottom row we also show the dimension $n_\eta$ of the radical.
\renewcommand{\arraystretch}{2}
$$
\begin{array}{|c||c|c|c|}
\hline
\eta & -\frac{1}{2m-2} & -\frac{1}{m-2}\\
\hline
n_\eta & 1 & m-1\\
\hline
\end{array}
$$
\end{theorem}

\pf We listed the critical values of $\eta$ and the corresponding radical dimensions in $M$ in the
beginning of this section. For $\eta=-\frac{1}{2m-2}$ (corresponding to the eigenvalue $k=2(n-2)$),
the dimension of the radical of $M$ is $1$, while $m_{2(n-2)}=0$. Hence $n_{-\frac{1}{2m-2}}=1-0=1$.
For $\eta=-\frac{1}{m-2}$ (corresponding to the eigenvalue $r=n-4$), the eigenspace dimension in $M$ is
$n-1=2m-1$, while $m_{n-4}=m_{2m-4}=m$. Hence the dimension of the radical of $Q_m(\eta)$
in this case is $2m-1-m=m-1$, as claimed. Finally, as we have already mentioned, $\eta=1$ (corresponding to the eigenvalue
$s=-2$) is excluded since $\eta\notin\{1,0\}$.\qed

\medskip
The second example of calculation will involve a flip of $M=M_\eta(G,D)$ for $G$ the symplectic group $Sp_{2m}(2)$. All flips for this $M$ have been determined in Joshi's PhD dissertation \cite{J} (and this will appear in \cite{JSS}). Here we determine the critical values and radical dimension for $M$ itself and for one series of flip subalgebras of $M$.

As in the symmetric case, the information for $M$ in the symplectic case is readily available from \cite{HS}. Namely, in the table in this paper, we find that the eigenvalues for the collinearity matrix $T$ are $2^{2m-1}$, $2^{m-1}$, and $-2^{m-1}$. Their multiplicities are
$1$, $2^{2m-1}-2^{m-1}-1$, and $2^{2m-1}+2^{m-1}-1$, respectively. Correspondingly, the critical values of $\eta$ for $M$ are $-\frac{2}{2^{2m-1}}=-\frac{1}{2^{2m-2}}$, $-\frac{2}{2^{m-1}}=-\frac{1}{2^{m-2}}$, and $-\frac{2}{-2^{m-1}}=\frac{1}{2^{m-2}}$, with the multiplicity being in each case the dimension of the radical. (As in the symmetric group case,
we may need to drop some eigenvalues for small $m$ since $\eta\notin\{1,0\}$.)

One interesting family of flip subalgebras in $M$ arises when the flip is induced by the rank $2$ symplectic element $\sg$ (this is $\sg_1$ in \cite{JSS}). In this case, the flip subalgebra $A(\sg)=M_\sg$ has a basis consisting of all of its $2^{2m-2}-1$ singles and all its $2^{2m-1}-2^{2m-3}$ doubles. In particular, the number $e$ of extras is zero, which also means that $M$ is the ambient algebra for $A(\sg)$. It follows that $\dim A_\sg=2^{2m-2}-1+2^{2m-1}-2^{2m-3}=2^{2m-1}+2^{2m-3}-1$ and $\dim W=\dim M-\dim A_\sg=2^{2m}-1-2^{2m-1}-2^{2m-3}+1=2^{2m-1}-2^{2m-3}=3\cdot 2^{2m-3}$. Therefore, our equations (\ref{eq 1}) and (\ref{eq 2}) become:
\begin{equation} \label{sum}
m_{2^{m-1}}+m_{-2^{m-1}}=3\cdot 2^{2m-3}
\end{equation}
and
\begin{equation} \label{difference}
2^{m-1}m_{2^{m-1}}+(-2^{m-1})m_{-2^{m-1}}=0.
\end{equation}
Equation (\ref{difference}) simplifies to $m_{2^{m-1}}-m_{-2^{m-1}}=0$, that is, $m_{2^{m-1}}=m_{-2^{m-1}}$. Comparing with (\ref{sum}), we get that $m_{2^{m-1}}=m_{-2^{m-1}}=3\cdot 2^{2m-4}$.

\begin{theorem}
Suppose that $M=M_\eta(G,D)$, where $G=Sp_{2m}(2)$, $m\geq 2$, and suppose that $\sg\in G$ is a symplectic element of rank $2$ and $A=M_\sg$. Then the radical of $A$ is non-zero if and only if $\eta$ is in the following table, where the second row shows the dimension of the radical:
\renewcommand{\arraystretch}{2}
$$
\begin{array}{|c||c|c|c|}
\hline
\eta & -\frac{1}{2^{2m-2}} & -\frac{1}{2^{m-2}} & \frac{1}{2^{m-2}}\\
\hline
n_\eta & 1 & 2^{2m-2}+2^{2m-4}-2^{m-1}-1 & 2^{2m-2}+2^{2m-4}+2^{m-1}-1\\
\hline
\end{array}
$$
\end{theorem}

\pf It is quite similar to the proof of Theorem \ref{orthogonal example}. The critical values are $-\frac{2}{k}=
-\frac{2}{2^{2m-1}}=-\frac{1}{2^{2m-2}}$, $-2\frac{2}{r}=-\frac{2}{2^{m-1}}=-\frac{1}{2^{m-2}}$, and $-\frac{2}{s}=
-\frac{2}{-2^{m-1}}=\frac{2}{2^{m-2}}$, as in the table. For the dimension of the radical, we get $1$ in the first case
(as always is the case for the eigenvalue $k$ of the adjacency matrix $T$), $(2^{2m-1}-2^{m-1}-1)-3\cdot 2^{2m-4}=
2^{2m-2}+2^{2m-4}-2^{m-1}-1$ for the eigenvalue $2^{m-1}$ and $(2^{2m-1}+2^{m-1}-1)-3\cdot 2^{2m-4}=
2^{2m-2}+2^{2m-4}+2^{m-1}-1$ for the eigenvalue $-2^{m-1}$ of $T$. This gives us the bottom row of the table.\qed

\medskip
When $m=2$, we have that $\frac{1}{2^{m-2}}=\frac{1}{1}=1$, and so this critical value needs to be omitted.

\section{Eigenvalues $0$ and $-1$}

These arise when the ambient Matsuo algebra is for a reducible $3$-transposition group. (Recall the discussion of irreducible vs reducible $3$-transposition groups in Section \ref{3-transposition}.) Let
us first consider a group $(G,D)$ with $\tau(G)\neq 1$, that is, $G=2^{\bullet h}\bar G$ for some $h\geq 1$ and
the smaller $3$-transposition group $(\bar G,\bar D)$, where $\bar G=G/\tau(G)$.

Consider $M=M_\eta(G,D)$ and define the subspace $U\subseteq M$ as follows:
$$U=\la c-d\mid c,d\in D,c\tau d\ra.$$

Recall the linear map $t:M\to M$ given by $t(d)=\sum_{c\in A_d}c$.

\begin{proposition}
The subspace $U$ is the $0$-eigenspace (i.e., the kernel) of $t$.
\end{proposition}

\pf If $c\tau d$ then $t(c-d)=t(c)-t(d)=\sum_{e\in A_c}e-\sum_{e\in A_d}e=0$ since $A_c=A_d$.
So $U$ is contained in the kernel of $t$. Conversely, we note that $M/U$ is naturally
isomorphic to $\bar M=M_\eta(\bar G,\bar D)$, where both $M$ and $\bar M$ are viewed simply as
vector spaces. Furthermore, if we define $\bar t$ to be the similar linear map on $\bar M$,
sending $\bar d\in\bar D$ to $\bar t(\bar d)=\sum_{\bar c\in A_{\bar d}}\bar c$, then the
action of $t$ on $\bar M\cong M/U$ coincides with $2^h\bar t$. So the claim now follows from
the fact that $\bar t$ does not have eigenvalue $0$ since $\tau(\bar G)=1$.\qed

\medskip
This has the following corollary, already noted in \cite{HS}.

\begin{corollary} \label{main tau}
Suppose that $G=2^{\bullet h}\bar G$ with $h\neq 0$ and let $\bar M$, $t$, $\bar t$ be as above. Then
the following hold:
\begin{enumerate}
\item[(a)] The non-zero eigenvalues of $t$ are of the form $\rho=2^h\bar\rho$, where $\bar\rho$ is
an eigenvalue of $\bar t$. Furthermore, the multiplicity of $\rho$ in $M$ is the same as the multiplicity
of $\bar\rho$ in $\bar M$.
\item[(b)] The multiplicity of $0$ in $M$ is $\dim(U)=\dim(M)-\dim(\bar M)=|D|-|\bar D|
=(2^h-1)|\bar D|$.
\end{enumerate}
\end{corollary}

A similar construction works in the case where $\theta(G)\neq 1$, that is, $G=3^{\bullet h}\hat G$,
where $h\geq 1$ and $(\hat G,\hat D)$ is a smaller $3$-transposition group, with $\hat G=G/\theta(G)$.
Similarly, we define the subspace $U=\la c-d\mid c,d\in D,c\theta d\ra$.

\begin{proposition} \label{-1-eigenspace}
$U$ is the $-1$-eigenspace of the map $t$.
\end{proposition}

\pf If $c\theta d$ then  $D_c=D_d$ and hence $A_c\cup\{c\}=A_d\cup\{d\}$. Thus,
$t(c-d)=t(c)-t(d)=\sum_{e\in A_c}e-\sum_{e\in A_d}e=-(c-d)$. So indeed, $U$ is contained in the
$-1$-eigenspace of $t$. Conversely, $M/U$ is naturally isomorphic to $\hat M=M_\eta(\hat
G,\hat D)$, where again we view both $M$ and $\hat M$ as vector spaces. Moreover, the action
of $t$ on $\hat M\cong M/U$ coincides with the action of $3^h\hat t+(3^h-1)Id$, where $\hat t$
is a similar map on $\hat M$ defined in terms of collinearity on $\hat D$. This is because
$A_d\cup\{d\}$ is the full preimage under the map $M\to\hat M$ of $A_{\hat d}\cup\{\hat d\}$.

It follows that the eigenvalues of $t$ on $\hat M$ are of the form $\rho=3^h\hat\rho+3^h-1$ for the
eigenvalues $\hat\rho$ of $\hat t$. Note that $3^h\hat\rho+3^h-1=-1$ if and only if $\hat\rho=-1$.
However, $-1$ is not an eigenvalue of $\hat t$ since $\theta(\hat G)=1$.\qed

\medskip
Again, we have the corollary that allows us to identify the multiplicities of eigenvalues.

\begin{corollary} \label{main theta}
Suppose that $G=3^{\bullet h}\hat G$ with $h\neq 0$ and let $\hat M$, $t$, $\hat t$ be as above. Then
the following hold:
\begin{enumerate}
\item[(a)] The eigenvalues of $t$, not equal to $-1$, are of the form $\rho=3^h\hat\rho+3^h-1$, where
$\hat\rho$ is an eigenvalue of $\hat t$. Furthermore, the multiplicity of $\rho$ in $M$ is the same
as the multiplicity of $\hat\rho$ in $\hat M$.
\item[(b)] The multiplicity of $-1$ in $M$ is $\dim(U)=\dim(M)-\dim(\hat M)=|D|-|\hat D|
=(3^h-1)|\hat D|$.
\end{enumerate}
\end{corollary}

These statements solve the eigenvalue and multiplicity problem for reducible Matsuo algebras and hence we can find all critical values and radical dimensions for them. Note that the eigenvalue $\rho=0$ of $T$, arising when $\tau(G)\neq 1$, does not lead to a critical value $\eta=-\frac{2}{\rho}$ of $M$ because we cannot divide by zero. On the contrary, the eigenvalue $\rho=-1$ of $T$, arising when $\theta(G)\neq 1$, does always lead to a critical value $\eta=-\frac{2}{-1}=2$ of $M$.

Let us now discuss what happens when we look at a flip subalgebra. Hence, let $\sg$ be
a flip of $M=M_\eta(G,D)$ for a reducible $3$-transposition group $(G,D)$, and let $A=A(\sg)$ be
the corresponding flip subalgebra. We assume that $M$ is ambient for $A$ and, furthermore,
$A=M_\sg$.

Let us first consider the case where $\tau(G)\neq 1$. Hence $G=2^{\bullet h}\bar G$ where $h\geq 1$ and $\bar G=G/\tau(G)$.
Recall the definition of $U=\la c-d\mid c,d\in D, c\tau d\ra$. We denote by $\pi$ the natural projection from $M$ onto $\bar M=M(\bar G,\bar D)$, where $\bar G=G/\tau(G)$ and $\bar D$ is the image of $D$ in $\bar G$. Recall that $\pi$ is just a linear map (and not an algebra homomorphism) and $U=\ker(\pi)$. Let $\bar t$ be the collinearity map on $\bar M$. We have already used above the following fact.

\begin{proposition} \label{almost commute}
We have that $\pi\circ t=2^h\,\bar t\circ\pi$.
\end{proposition}

Similarly, we also have a statement for the flip $\sg$.

\begin{proposition} \label{they commute}
There is a flip $\bar\sg$ of $\bar M$ such that $\pi\circ\sg=\bar\sg\circ\pi$.
\end{proposition}

\pf Note that $\sg$, being an automorphism of the Fischer space of $(G,D)$, preserves collinearity and hence takes $\tau$-equivalence classes to $\tau$-equivalence classes. Since these equivalence classes naturally correspond to the points of the Fischer space of $(\bar G,\bar D)$, we have a permutation $\bar\sg$ on $\bar D$, such that $\pi\circ\sg=\bar\sg\circ\pi$ holds. Clearly, $\bar\sg$ preserves collinearity on $\bar D$ and the order of $\bar\sg$ divides $|\sg|=2$, so $\bar\sg$ is a flip or, possibly, the identity map.
\qed

\medskip
Turning to the flip subalgebra, we can claim the following.

\begin{proposition} \label{good map}
$\pi(M_\sg)=\bar M_{\bar\sg}$ and $\pi([M,\sg])=[\bar M,\bar\sg]$.
\end{proposition}

\pf Clearly, if $u\in M_\sg$ (respectively, $u\in [M,\sg]$) then $\bar u\in\bar M_{\bar\sg}$ (respectively, $\bar u\in[\bar M,
\bar\sg]$). Indeed, if $u^\sg=u$ then $\bar u^{\bar\sg}=\pi(u)^{\bar\sg}=\pi(u^\sg)=\pi(u)=\bar u$, where we used
Proposition \ref{they commute}. Similarly, if $u=[v,\sg]=v-v^\sg$ for some $v\in M$ then $\pi(u)=\pi(v-v^{\sg})=\pi(v)-\pi(v^\sg)=\pi(v)-\pi(v)^{\bar\sg}=[\bar v,\bar\sg]$. So $\pi$ maps $M_\sg$ to $\bar M_{\bar\sg}$ and $[M,\sg]$ to $[\bar M,\bar\sg]$.

It remains to show that both maps are onto. Consider a basis element $\bar v\in\bar M_{\bar\sg}$, i.e., $\bar v$ is a single, double, or extra in $\bar M$. If $\bar v=\bar c$ is a single, for some $\bar c\in\bar D$, then select $c\in D$ such that $\bar v=\bar c=\pi(c)$. If $c$ is itself a single in $M$ then, clearly, $c=c^\sg\in M_\sg$ and so $\bar v=\pi(c)\in\pi(M_\sg)$. If $c$ is not a single then $c+c^\sg$ is a double, and so $\bar v=\frac{1}{2}(\bar v+\bar v^{\bar\sg})=\frac{1}{2}\pi(c+c^\sg)\in\pi(M_\sg)$, where we used that $\bar v=\bar v^{\bar\sg}$ since $\bar v$ is a single. If $\bar v=\bar c+\bar c^{\bar\sg}$ is a double or extra, with $\bar c\in\bar D$, we again select $c\in D$ such that $\pi(c)=\bar c$ and then $c+c^\sg$ is a double or an extra and $\pi(c+c^\sg)=\pi(c)+\pi(c^\sg)=\pi(c)+\pi(c)^{\bar\sg}=\bar v$. Hence $\bar v$ is in the image under $\pi$ of $M_\sg$.
A very similar argument applies to show that $\pi$ maps $[M,\sg]$ onto $[\bar M,\bar\sg]$. Here instead of singles, doubles, and extras, we use skew doubles and skew extras.\qed

\medskip
One could ask what $\pi$ does to singles, doubles and extras. As we have already seen in the above proof, $\pi$ maps singles to singles. More generally, it maps $\Sg$-orbits to $\bar\Sg$-orbits, where $\bar\Sg=\la\bar\sg\ra$. In particular, $\pi$ maps a double $a+a^\sg$ to the double $\bar a+\bar a^{\bar\sg}$ if $\bar a\neq\bar a^{\bar\sg}$ and, otherwise, $\pi(a+a^\sg)=2\bar a$. Here, of course, $\bar a=\pi(a)$. For an extra $a+a^\sg$, we only have one option: $\pi(a+a^\sg)=\bar a+\bar a^{\bar\sg}$ is necessarily an extra in $\bar M$, since we cannot have in this case that $\bar a=\bar a^{\bar\sg}$. Indeed, if $\pi(a)=\pi(a^\sg)$ then $a\tau a^\sg$ and so $a$ and $a^\sg$ are non-collinear.

Now we aim to prove a statement about flip subalgebras which is roughly equivalent to Corollary \ref{main tau}.

\begin{proposition} \label{eigenspace map}
Suppose that $G=2^{\bullet h}\bar G$ with $h\neq 0$ and let $\bar M$, $t$, $\bar t$, $\sg$, $\bar\sg$, $\pi$ and $U$ be as above. Then the following hold:
\begin{enumerate}
\item[(a)] if $\rho\neq 0$ is an eigenvalue of $t$ then $\pi$ maps the $\rho$-eigenspace of $t$ in $M_\sg$ isomorphically onto the $\frac{\rho}{2^h}$-eigenspace of $\bar t$ in $\bar M_{\bar\sg}$;
\item[(b)] the $0$-eigenspace of $t$ in $M_\sg$ coincides with $M_\sg\cap U$ and it has dimension $\dim(M_\sg)-\dim(\bar M_{\bar\sg})$.
\end{enumerate}
\end{proposition}

\pf Indeed, since the $0$-eigenspace of $t$ is $U$, we immediately obtain the first claim in (b). Furthermore, since $U$ is the
kernel of $\pi$ and since $\pi(M_\sg)=\bar M_{\bar\sg}$ by Proposition \ref{good map}, we conclude that $\dim(\bar M_{\bar\sg})=\dim(M_\sg)-\dim(\ker(t|_{M_\sg}))=\dim(M_\sg)-\dim(M_\sg\cap U)$. This proves (b).

For (a), note that if $u\in M_\sg$ is a $\rho$-eigenvector of $t$ then $\bar t(\pi(u))=\frac{1}{2^h}\pi(t(u))=\frac{1}{2^h}\pi(\rho u)=\frac{\rho}{2^h}\pi(u)$. (We used Proposition \ref{almost commute}.) Hence $\pi$ maps the $\rho$-eigenspace of $t$ to the $\frac{\rho}{2^h}$-eigenspace of $\bar t$.  It remains to see that this map is an isomorphism. First of all, this restriction map is injective since $\ker(\pi)=U$, being the $0$-eigenspace of $t$, intersects trivially with the $\rho$-eigenspace of $t$. Furthermore, note that $t$ is semisimple on both $M$ and $M_\sg$ and, similarly, $\bar t$ is semisimple on $\bar M$ and $\bar M_{\bar\sg}$. Taking a $\frac{\rho}{2^h}$-eigenvector $\bar u\in\bar M_{\bar\sg}$ of $\bar t$, let $u\in M_\sg$ be such that $\pi(u)=\bar u$. By semisimplicity, $u=u_0+u_{\rho_1}+\ldots+u_{\rho_k}$, where $u_0$ is a $0$-eigenvector of $t$ in $M_\sg$ and, similarly, every $u_{\rho_i}$ is a $\rho_i$-eigenvector of $t$ in $M_\sg$. Then $\bar u=\pi(u)=\pi(u_0)+\pi(u_{\rho_1})+\ldots+\pi(u_{\rho_k})$. By the above, $\pi(u_0)=0$ and each $\pi(u_{\rho_i})$ is a $\frac{\rho_i}{2^h}$-eigenvector of $\bar t$. Since $\bar u$ has been chosen to be itself an eigenvector, we conclude that $\bar u$ coincides with $\pi(u_{\rho_i})$, where $\rho_i=\rho$, and all other summands are zero. In particular, this means that $\pi$ maps the $\rho$-eigenspace of $t$ in $M_\sg$ onto the $\frac{\rho}{2^h}$-eigenspace of $\bar t$ in $\bar M_{\bar\sg}$. Hence (a) holds.\qed

\medskip
Similar ideas apply in the case of groups $G=3^{\bullet h}\hat G$. That is, here $\theta(G)\neq 1$ and $\hat G=G/\theta(G)$. Let again $\pi$ be the natural surjective linear map from $M=M_\eta(G,D)$ to $\hat M=M_\eta(\hat G,\hat D)$ and let $U=\ker(\pi)$. Then, as we discussed earlier, $U$ is the
$-1$-eigenspace of the collinearity map $t$ and $\dim(U)=|D|-|\hat D|=(3^h-1)|\hat D|$. Let $\hat t$ be the collinearity map on $\hat M$.
The following fact has already been mentioned in the proof of Proposition \ref{-1-eigenspace}. (This is an analogue of Proposition \ref{almost commute}.)

\begin{proposition} \label{almost commute -1}
We have that $\pi\circ t=(3^h\hat t+(3^h-1)Id_{\hat M})\circ\pi$.
\end{proposition}

Similarly, we have an analogue of Proposition \ref{they commute}.

\begin{proposition} \label{they commute -1}
There is a flip $\hat\sg$ of $\hat M$ such that $\pi\circ\sg=\hat\sg\circ\pi$.
\end{proposition}

The proof is essentially identical to that of Proposition \ref{they commute}, with $\tau$-equivalence classes replaced with $\theta$-equivalence classes, and hence also bars replaced with hats.

\begin{proposition} \label{good map -1}
$\pi(M_\sg)=\hat M_{\hat\sg}$ and $\pi([M,\sg])=[\hat M,\hat\sg]$.
\end{proposition}

Again, the proof of this is essentially identical to that of Proposition \ref{good map}, but note that in the case, where $\hat v=\hat c$ is a single, for some $\hat c\in\hat D$, the point $c\in D$ satisfying $\pi(c)=\hat c$ is either a single or part of an extra, not a double, as in the proof of Proposition \ref{good map}.

More in particular, $\pi$ maps singles from $M_\sg$ to singles from $\hat M_{\hat\sg}$ and it maps doubles from $M_\sg$ to doubles from $\hat M_{\hat\sg}$. For an extra $c+c^\sg\in M_\sg$, we have that $\pi(c+c^\sg)=\pi(c)+\pi(c)^{\hat\sg}$ is an extra when $\pi(c)\neq\pi(c)^{\hat\sg}$, and $\pi(c+c^\sg)=2\pi(c)$ when $\pi(c)=\pi(c)^{\hat\sg}$ (i.e., $\pi(c)$ is a single in this case).

We are now prepared for the main statement.

\begin{proposition} \label{eigenspace map -1}
Suppose that $G=3^{\bullet h}\hat G$ with $h\neq 0$ and let $\hat M$, $t$, $\hat t$, $\sg$, $\hat\sg$, $\pi$ and $U$ be as above. Then the following hold:
\begin{enumerate}
\item[(a)] if $\rho\neq -1$ is an eigenvalue of $t$ then $\pi$ maps the $\rho$-eigenspace of $t$ in $M_\sg$ isomorphically onto the $(\frac{\rho+1}{3^h}-1)$-eigenspace of $\hat t$ in $\hat M_{\hat\sg}$;
\item[(b)] the $-1$-eigenspace of $t$ in $M_\sg$ coincides with $M_\sg\cap U$ and it has dimension $\dim(M_\sg)-\dim(\hat M_{\hat\sg})$.
\end{enumerate}
\end{proposition}

\pf Part (b) is immediate since $U$ is at the same time the $-1$-eigenspace of $t$ in $M$ and the kernel $\ker(\pi)$.

Turning to (a), if $\rho\neq -1$ is an eigenvalue and $u\in M_\sg$ is a $\rho$-eigenvector of $t$ then $\hat t(\pi(u))=\frac{1}{3^h}(\pi(t(u))-(3^h-1)\pi(u))=\frac{1}{3^h}(\pi(\rho u))-(3^h-1)\pi(u))=\frac{\rho-3^h+1}{3^h}\pi(u)=(\frac{\rho+1}{3^h}-1)\pi(u)$, where we utilised Proposition \ref{almost commute -1}. Hence $\pi$ maps the $\rho$-eigenspace of $t$ to the
$\mu$-eigenspace of $\bar t$, where $\mu=\frac{\rho+1}{3^h}-1$. The proof of this restriction being a linear isomorphism
is virtually identical to the final part of the proof of Proposition \ref{eigenspace map}.\qed

\section{A symplectic example}

We can illustrate the case where we have eigenvalue $0$ using examples of flip subalgebras arising in the symplectic case.

Recall that the symplectic group $G=Sp_{2m}(2)$ acts on the non-degenerate symplectic space $V$ of dimension $n=2m$. It is a $3$-transposition group where the class $D$ of $3$-transpositions is the class of symplectic transvections $r_v:V\to V$, $v\in V^\sharp=V\setminus\{0\}$, where $r_v(u)=u+(u,v)v$ for all $u\in V$. Here $(\cdot,\cdot)$ is the symplectic form on $V$.

It will be convenient to identify each transvection $r_v$ with the corresponding vector $v\in V^\sharp$, so that the point set of the Fischer space (and the set of axes of the Matsuo algebra $M$) is $V^\sharp$. The lines of the Fischer space are triples $\{u,v,u+v\}$ such that $(u,v)\neq 0$. Since we have two additions for vectors from $V^\sharp$, the one in $V$ and the other one in $M$, we will use the different symbol $\dotplus$ for the latter addition.

The flips of the symplectic Matsuo algebra $M=M(G,D)$ have been classified in \cite{JSS}. Since we will need to identify the flip $\bar\sg$ for
the selected $\sg$, we will need some details of this classification. Each flip $\sg$ of $M_\eta(Sp_{2m}(2))$ is identified, up to
conjugation, by its type and rank. The \emph{rank} of $\sg$ is simply the dimension of the commutator subspace $[V,\sg]$.
The \emph{type} of $\sg$ is defined in terms of the so-called second form $\llf u,v\rrf=(u,v^\sg)$, which is a symmetric bilinear form on $V$. We also define $V(\sg)=\{v\in V\mid\llf v,v\rrf=0\}$, a subspace of $V$ which is either a hyperplane (flip type $1$) or
the entire $V$ (flip type $2$). For flips of type $1$, the rank ranges from $1$ to $m$, while for type $2$, the rank can only be
even and again between $1$ and $m$. Since the induced flip $\bar\sg$ can sometimes be the identity, we consider the latter as
an ``honorary'' flip, which is then of type $2$ and rank $0$.

Representatives of all conjugacy classes of flips can be constructed as follows. Let $E=\la e_1,e_2,\ldots,e_m\ra$ be a maximal totally isotropic subspace of $V$. Then $\tau_i=r_1r_2\cdots r_i$, where $r_j=r_{e_j}$, is of type $1$ and rank $i$. Representatives of classes of flips of type $2$ are $\sg_i=s_1s_2\ldots s_i$, where $s_j=r_{e_{2i-1}}r_{e_{2i}}r_{e_{2i-1}+e_{2i}}$. The flip $\sg_i$ has rank $2i$.

Now we are fully prepared to give an example involving the eigenvalue $0$. It arises in the flip subalgebra when the ambient Matsuo subalgebra $M_0$ is a proper subalgebra of $M$. This happens in the symplectic case when $\sg$ is of type $1$. E.g., we can take, say, $\sg=\tau_m=r_1r_2\cdots r_m$. According to \cite{JSS}, when $\sg=\tau_m$, all singles and doubles (and no extras) arise from the hyperplane $V_0=V(\sg)=e^\perp$, where $e=e_1+e_2+\ldots+e_m$. Hence $M_0$ is the Matsuo subalgebra of $M$ corresponding to the subspace $V_0^\sharp$ of the Fischer space $V^\sharp$. Note that $V_0^\perp=\la e\ra$ is the $1$-dimensional radical of $V_0$ and hence the ambient algebra $M_0$ is the direct sum of the $1$-dimensional algebra $\lla e\rra$ and the Matsuo algebra $M_0':=\lla V_0\setminus\la e\ra\rra$.

Clearly, $\sg$ leaves both $e$ and $V_0\setminus\la e\ra$ invariant and so the flip subalgebra $A=A(\sg)$ also decomposes as the direct sum of $\lla e\rra$ and $A'=A(\sg)'=(M_0')_\sg$. The $1$-dimensional summand is, clearly, simple, so the radical of $A$ coincides with the radical of $A'$. Hence we just investigate the latter. For this subalgebra, its ambient algebra is $M_0'$, and it is the Matsuo algebra of the group $G_0'\cong 2^{2m-2}:Sp_{2m-2}(2)$, which arises as a maximal parabolic subgroup of $G=Sp_{2m}(2)$. Then $G_0'=2^{\bullet 1}:Sp_{2m-2}(2)$ in the notation of \cite{HS}. According to the second big table in that paper, the eigenvalues of the collinearity map in this case are $2^{2m-2}$ (multiplicity $1$), $2^{m-1}$ (multiplicity $2^{2m-3}-2^{m-2}-1$), $-2^{m-1}$ (multiplicity $2^{2m-3}+2^{m-2}-1$), and $0$ (multiplicity $2^{2m-2}-1$). This gives us the critical values $-\frac{2}{2^{2m-2}}=-\frac{1}{2^{2m-3}}$, $-\frac{2}{2^{m-1}}=-\frac{1}{2^{m-2}}$, and $-\frac{2}{-2^{m-1}}=\frac{1}{2^{m-2}}$. Note that the eigenvalue $0$ does not lead to a critical value.

According to the method we developed in the preceding section, we need to determine the induced flip $\bar\sg$ of $\bar M_0=M_\eta(\bar G_0,\bar D_0)$, where $\bar G_0=G_0/\tau(G_0)\cong Sp_{2m-2}(2)$.

\begin{proposition}
The induced flip $\bar\sg$ is of type $2$ and rank $m-1$, if $m$ is odd, and $m-2$, if $m$ is even.
\end{proposition}

\pf The symplectic space underlying $\bar G_0$ is $V_0/\la e\ra$. Take an arbitrary $u\in V_0$. Then $\llf\bar u,\bar u\rrf=(u,u^\sg)=(u,u+\sum_{j=1}^m(u,e_j)e_j)=(u,u+\sum_{j=1}^me_j)=(u,u+e)=(u,u)=0$. Here we used the fact that the terms we added to the sum gave the form value zero. We also used, for the last equality, that the form is symplectic. We have shown that $\llf\bar u,\bar u\rrf=0$ for all $\bar u\in\bar V_0$, which means that $V(\bar\sg)=\bar V_0$, i.e., $\bar\sg$ is of type $2$.

It remains to determine the rank of $\bar\sg$. Clearly, $[\bar V_0,\bar\sg]=\overline{[V_0,\sg]}$. Note that $[V,\sg]=E$ is of dimension $m$. Furthermore, the linear map $V\to [V,\sg]=E$ given by $u\mapsto [u,\sg]=u^\sg-u$ has kernel $E$. Since $V_0$ is a hyperplane of $V$ and $V_0\supseteq E$, the rank-plus-nullity theorem gives us that the dimension of $[V_0,\sg]$ is $m-1$. It is now easy to see from the exact formula of this linear map that $[V_0,\sg]$ consists of all $u\in E$ having even support within the basis $\{e_1,e_2,\ldots,e_m\}$. When $m$ is odd, $e\notin[V_0,\sg]$ and so $\dim[\bar V_0,\bar\sg]=\dim[V_0,\sg]=m-1$. On the other hand, when $m$ is even, $e\in[V_0,\sg]$ and so $\dim[\bar V_0,\bar\sg]=\dim[V_0,\sg]-1=m-2$.\qed

\medskip
This means that our $\bar\sg$ is conjugate to $\sg_j$ with
$j=\lfloor\frac{m-1}{2}\rfloor$, this formula covering both cases above. Now that
we identified $\bar\sg$, we can find the multiplicities $m_{\bar\theta_i}$ of the
eigenvalues $\bar\theta_i=\frac{1}{2}2^{2m-2}=2^{2m-3}$,
$\frac{1}{2}2^{m-1}=2^{m-2}$, and $\frac{1}{2}(-2^{m-1})=-2^{m-2}$ of the
collinearity map $\bar t$ within $\bar W=[\bar M_0,\bar\sg]$. (Note that the values of $\bar\theta_i$ are in agreement with the tables in \cite{HS}.)

First of all, we will need the following information about $\sg_i$ from \cite{JSS}.
For $\sg_i$, there are no extras within $\bar D_0$. Furthermore, there are
$2^{2m-2j-2}-1$ singles and $2^{2m-3}-2^{2m-2j-3}$ doubles. In particular,
$$
\dim\bar W=\left\{\begin{array}{rl}2^{2m-3}-2^{m-2},&\mbox{if $m$ is odd};\\
2^{2m-3}-2^{m-1},&\mbox{if $m$ is even};\end{array}\right.
$$
since $\dim\bar W$ is the sum of numbers of doubles and extras and since $2j$ is
the rank of $\bar\sg$.

Next, we have the following equations tying the multiplicities $m_{\bar\theta_i}$
together. Equations (\ref{eq 1}) and (\ref{eq 2}) give us:
\begin{equation}\label{dim}
m_{2^{m-2}}+m_{-2^{m-2}}=
\left\{\begin{array}{rl}2^{2m-3}-2^{m-2},&\mbox{if $m$ is
odd};\\2^{2m-3}-2^{m-1},&\mbox{if $m$ is even}; \end{array}\right.
\end{equation}
and
\begin{equation}\label{weighted}
m_{2^{m-2}}2^{m-1}-m_{-2^{m-2}}2^{m-2}=m_{2^{m-2}}2^{m-2}-m_{-2^{m-2}}2^{m-2}=0
\end{equation}
(Indeed, the right side of the first equation is $\dim\bar W$ and for the second
equation it is zero, since there are no extras within $\bar M_0$.) The second
equation simplifies to $m_{2^{m-2}}-m_{-2^{m-2}}=0$, that is,
$m_{2^{m-2}}=m_{-2^{m-2}}$. Now from the first equation we get that
$$
m_{2^{m-2}}=m_{-2^{m-2}}=\left\{\begin{array}{rl}2^{2m-4}-2^{m-3},&\mbox{if $m$
is odd};\\2^{2m-4}-2^{m-2},&\mbox{if $m$ is even}.\end{array}\right.
$$
(Note that $m\geq 2$ and so these formulas make sense.)

We can now formulate the following.

\begin{theorem}
Suppose that $M=M_\eta(G,D)$, where $G=Sp_{2m}(2)$, $m\geq 2$, and suppose that the flip $\sg=\tau_m$ is of type $1$ and rank $m$. Let $A=A(\sg)$. Then the radical of $A$ is non-zero if and only if $\eta$ is one of the critical values in the table below, where the second row shows the dimension of the radical:
\renewcommand{\arraystretch}{2}
$$
\begin{array}{|c||c|c|c|}
\hline
\eta & -\frac{1}{2^{2m-3}} & -\frac{1}{2^{m-2}} & \frac{1}{2^{m-2}}\\
\hline
n_\eta & 1 & 2^{2m-4}-2^{m-3}-1,\mbox{ $m$ odd,}& 2^{2m-4}+3\cdot 2^{m-3}-1,\mbox{ $m$ odd},\\
         &    & 2^{2m-4}-1,\mbox{ $m$ even}              & 2^{2m-4}+2^{m-1}-1,\mbox{ $m$ even}\\
\hline
\end{array}
$$
\end{theorem}

\pf The critical values were discussed earlier, so the only thing to verify is the dimensions of the corresponding eigenspaces of $t$ within $A$. These are the multiplicities of the eigenvalues $\theta$ within $M_0$ minus the corresponding $m_\theta=m_{\bar\theta}$. The calculation yields the values in the table. Note again that the eigenvalue $\theta=0$ does not lead to a critical value.\qed

\medskip
Unfortunately, we do not have an example illustrating application of our method when we have the additional eigenvalue $-1$. This is because the case of orthogonal groups in characteristic $3$ has not been looked at yet. However, this case follows essentially the same scheme as above.

\section{Positive characteristic}

In the above results we rarely mentioned the characteristic of the ground field $\F$. However, the characteristic does matter. The results from \cite{HS}, that we used in our exams are obtained for the case of characteristic zero. The purpose of this brief section is to discuss whether the same results hold in at least some positive characteristics.

Note that the Matsuo algebra $M=M_\eta(G,D)$ is defined as long as the characteristic of $\F$ is not $2$. Furthermore, the version of $M$ is positive characteristic $p\neq 2$ can be obtained by reducing modulo $p$ the suitable characteristic $0$ version of $M$. If we want our method to work in characteristic $p$, the following key fact must remain true: the collinearity map $t$ is semisimple,
i.e., it admits a basis of eigenvectors. Let $\theta_1,\theta_2,\ldots,\theta_r$ are the eigenvectors of $t$ in characteristic $0$. (Note that they are integer and so they can be reduced modulo any $p$.) Since $t$ is semisimple, its minimal polynomial is integral
and it has no multiple roots. Let us call $p$ \emph{good} if no two eigenvalues $\theta_i$ and $\theta_j$, $i\neq j$, are congruent modulo $p$.

\begin{proposition}
If the odd prime $p$ is good, then the collinearity map $t$ remains semi-simple in characteristic $p$.
\end{proposition}

\pf Since the characteristic $p$ version of $M$ is obtained by reducing the characteristic $0$ version modulo $p$, the minimal polynomial $f$ of $t$ (for characteristic $0$) reduced modulo $p$ annihilates the mod-$p$ version of $t$. But if $p$ is good,
then the roots of $f$ remain distinct modulo $p$, i.e. it remains a multiplicity-free polynomial. So the minimal polynomial of the mod-$p$ version of $t$ is also multiplicity-free, i.e., $t$ remains semi-simple.\qed

\medskip
We now claim that the remaining part of the calculation is not sensitive to the characteristic of $\F$. Indeed, all ingredients of our calculation, such as the multiplicities of eigenvalues, numbers of singles, doubles and extras, the linear equations (\ref{eq 1}), (\ref{eq 2}) and (\ref{eq 3}), are not sensitive to the field characteristic, as it simply never came up. The conclusion to this is
the following main result.

\begin{theorem}
Suppose that the ground field $\F$ is of odd characteristic $p$ that is good for the $3$-transposition group $(G,D)$. Then the critical values of $M=M_\eta(G,D)$ can be obtained by reducing modulo $p$ the critical values of the characteristic $0$ version of $M$ and the dimension of the radical in each case remains unchanged after reducing modulo $p$.
\end{theorem}

The ``bad'' characteristics $p$ are more difficult and need to be treated individually, but the good news is that there are only finitely many of them for each $(G,D)$.

\begin{example}
According to \cite{HS}, if $G=SU_n(4)$ then, in characteristic $0$, the eigenvalues of the collinearity map $t$ are: $2^{2n-3}$, $-(-2)^{n-3}$ and $-(-2)^{m-2}$. So the ``bad'' primes in this case are the odd primes $p$ that divide $2^{2n-3}+(-2)^{n-3}=2^{n-3}(2^n+(-1)^{n-3})=2^{n-2}(2^n-(-1)^n)$ or $2^{2n-3}+(-2)^{n-2}=2^{n-3}(2^n+(-1)^{n-2})=2^{n-2}(2^n+(-1)^n)$. So they are the primes dividing $2^n\pm 1$.
\end{example}


\begin{thebibliography}{99}
	
\bibitem{ABS} M.~Alsaeedi, V.~Bovdi, and S.~Shpectorov, Flip subalgebras of Matsuo algebras for extended symmetric groups, in preparation.

\bibitem{CH} H.~Cuypers and J.I.~Hall, The $3$-transposition groups with trivial center, {\it J. Algebra} {\bf 178} (1995), 149--193.

\bibitem{F} B.~Fischer, Finite groups generated by $3$-transpositions. I., {\it Invent. Math.} {\bf 13} (1971), 232--246.

\bibitem{DF} D. Fox, The commutative nonassociative algebra of metric curvature tensors, \textit{Forum Math. Sigma} {\bf 9} (2021), Paper No. e79, 48 pages.

\bibitem{GJMSS} A.~Galt, V.~Joshi, A.~Mamontov, S.~Shpectorov and A.~Staroletov, Double axes and subalgebras of Monster type in Matsuo algebras, \textit{Comm. Algebra} {\bf 49}, no. 10, 4208--4248.

\bibitem{GS} I.~Gorshkov and A.~Staroletov, private communication.

\bibitem{HRS1} J.I.~Hall, F.~Rehren, and S.~Shpectorov, Universal axial algebras and a theorem of Sakuma, \textit{J. Algebra} {\bf 421} (2015), 394--424.

\bibitem{HRS2} J.I.~Hall, F.~Rehren, and S.~Shpectorov, Primitive axial algebras of Jordan type, \textit{J. Algebra} {\bf 437} (2015), 79--115.

\bibitem{HSS2} J.I.~Hall, Y.~Segev, and S.~Shpectorov, On primitive axial algebras of Jordan type, \textit{Bull. Inst. Math. Acad. Sin. \textup{(}N.S.\textup{)}}, {\bf 13} (2018), no. 4, 397--409.

\bibitem{HS} J.I.~Hall and S.~Shpectorov, The spectra of finite $3$-transposition groups, {\it Arab. J. Math.}, {\bf 10} (2021), 611--638.

\bibitem{HRS3} C.G.~Hoffman, B.G.~Rodrigues, and S.~Shpectorov, Axial algebras of Monster type from unitary groups, in preparation.

\bibitem{J} V.~Joshi, {\it Axial Algebras of Monster type $(2\eta,\eta)$}, PhD thesis, University of Birmingham, 2020.

\bibitem{JSS} V.~Joshi, Y.~Shi and S.~Shpectorov, Flip subalgebras in Matsuo algebras of symplectic
and orthogonal type, in preparation.

\bibitem{KMS} S.M.S.~Khasraw, J.~McInroy, S.~Shpectorov, On the structure of axial algebras, {\it Trans. Amer. Math. Soc.} {\bf 373} (2020), 2135--2156.

\bibitem{M} A. Matsuo, {\it $3$-transposition  groups of symplectic type
and vertex operator algebras} (version 1), manuscript, November 2003
(available as {\tt arXiv.math/0311400v1}).

\bibitem{MS} J.~M\textsuperscript{c}Inroy and S.~Shpectorov, Axial algebras of Jordan and Monster type, {\it Proc. Groups
St Andrews $2022$}, to appear. https://arxiv.org/abs/2209.08043

\bibitem{T} V.G.~Tkachev, The universality of one half in commutative nonassociative algebras with identities, {\it J. Algebra} {\bf 569} (2021), 466--510.


\end{thebibliography}
\end{document}